\documentclass[11pt]{amsart}

\usepackage{url}
\input epsf

\newtheorem{Thm}{Theorem}
\newtheorem{Lem}{Lemma}

\newcommand{\whbox}{$\square$ }

\begin{document}
\thanks{Work on this material was done at home and at Wolfson College, Oxford. I am grateful for Oxford hospitality and Winnipeg patience. I also want to acknowledge the help, with my first attempt at studying this material, of Lorrita McKnight, then a graduate student of the Department of Applied Mathematics at the University of Manitoba, and Mark Fleming of IT Services at Queen's University, Kingston.}

\title[Multicoloured Isonemal Fabrics by Thick Striping]{Colouring Isonemal Fabrics with more than two Colours by Thick Striping}

\author{R.~S.~D.~Thomas}

\address{St John's College and Department of Mathematics, University of Manitoba, Winnipeg, Manitoba  R3T 2N2  Canada}

\email{thomas@cc.umanitoba.ca}

\begin{abstract}
Perfect colouring of isonemal fabrics by thin and thick striping of warp and weft with more than two colours is introduced. Conditions that prevent perfect colouring by striping are derived, and it is shown that avoiding them is sufficient to allow it. Examples of thick striping in all possible species are given.
\end{abstract}

\maketitle

\section{Introduction}%
\label{sect:Intro}

\noindent This paper intends to introduce the topic of weaving with more than two colours into the mathematical literature and to derive some initial general results. The first section introduces weaving and striping, the second discusses colouring without restriction on the number of colours, the third the restriction to the same finite number of colours in both directions, the fourth the example of thick striping with three colours, and the fifth the same with four and six colours.

All depends on the perfectness of the fabric colouring. What has to be arranged is that the colouring of the strands of the fabric be compatible with the symmetry group of the fabric, so that every symmetry operation of the underlying fabric's topology is a coherent mapping of the colours of the strands. The colours need not be preserved but may be permuted.

The first result (Lemma \ref{lem:1}) shows what restrictions are needed for thick striping with all stripes of different colours. No restrictions (Lemma \ref{lem:2}) allow perfect colouring by the assignment of different colours to all single strands (thin striping) with all stripes of different colours. All sorts of fabrics (Theorem \ref{thm:1}) can be perfectly colour\-ed by thin striping or, with the exception of fabrics in genus I or II, assignment of different colours to all adjacent pairs of strands (thick striping). Section \ref{sect:Finite} restricts the number of colours to finite and finds a couple of necessary conditions on the colour choice for perfect colouring, which jointly are sufficient (Theorem \ref{thm:2}). The colours of the warps must be either the same as the colours of the wefts or they must be all different (Lemma \ref{lem:5}). More detailed results state that in certain species of fabrics it is possible to colour warps and wefts perfectly by thick striping with an odd number of colours (Theorem \ref{thm:4}) and in certain species of fabrics it is possible with an even number of colours (Theorem \ref{thm:5}).
\section{Background}
\label{sect:Background}

\noindent Weaving was introduced into the mathematical literature by Gr\"unbaum and Shephard with `Satins and twills: An introduction to the geometry of fabrics' \cite{GS1980}, and satins and twills are a good place to begin.
First, a \emph{fabric} is a woven structure that is sufficiently coherent that it does not fall apart.
Without that constraint it is called a \emph{prefabric} \cite{GS1988}.
Basic definitions are in \cite{GS1980} or as indicated.
The weaving patterns of interest to mathematicians so far have been those that Gr\"unbaum and Shephard introduced, periodic in the plane with strands perpendicular and with symmetry groups $G_1$ transitive on strands, making them \emph{isonemal}.
The symmetry groups to which reference will be made here are the wallpaper groups $G_1$ of \cite{Roth1993} rather than the full three-dimensional symmetry group $G$ of the prefabric.
Each symmetry operation of $G_1$ will be combined with (or not) reflection in the plane $E$ of the fabric, represented by $\tau$.
In the standard diagram, called the \emph{design}, of a weaving pattern, the vertical strands called \emph{warps} are coloured dark and the horizontal strands called \emph{wefts} are coloured pale so that the visual appearance indicates which strand is over which strand in the places where they cross.
The places where they cross are all that is of interest, and so they are taken to be square \emph{cells} without boundary tessellating the plane, infinite for convenience.
So much is square here that the term `cell' was introduced in \cite{Thomas2009}.
The requirement of periodicity means that there are fundamental blocks of the pattern whose translates cover the plane \cite{Schatt1978}.
The requirement of transitivity of the symmetry group/isonemality of the fabric means that every strip of cells can be transformed into every other strip of cells, vertical or horizontal, by a symmetry operation, perhaps with the reversal of dark and pale, effected by $\tau$ \cite{Roth1993}.
The standard colouring of the strands of an isonemal weaving structure is a \emph{perfect colouring} in the standard sense that every symmetry operation preserves the colouring.
The definition of perfect colouring allows reversal of dark and pale; it needs to preserve the colouring rather than the colours.
The operation $\tau$ cannot itself be a symmetry operation, but it is involved in many symmetry operations; it represents a total reversal of which strand covers which strand over the whole plane, reflection in $E$, the plane of the fabric.
I illustrate these ideas with a few examples.
\begin{figure}
\[
\begin{array}{cc}
\epsffile{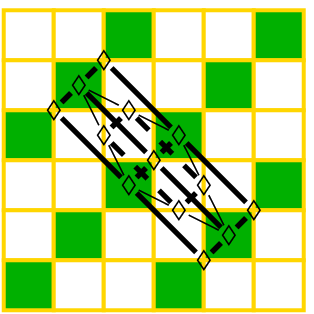} &\epsffile{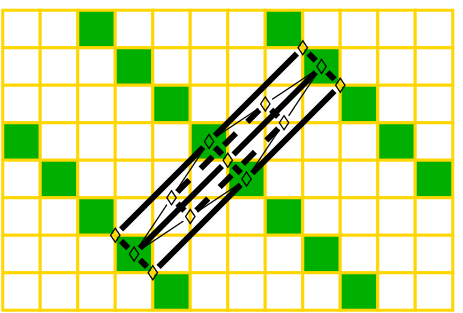}\\
\mbox{(a)} &\mbox{(b)}
\end{array}
\]
\caption{Examples of twills to be used later. a. 2/1 twill. b. 4/1 twill.}
\label{fig:8d;8e}
\end{figure}

\begin{figure}
\[
\begin{array}{cc}
\epsffile{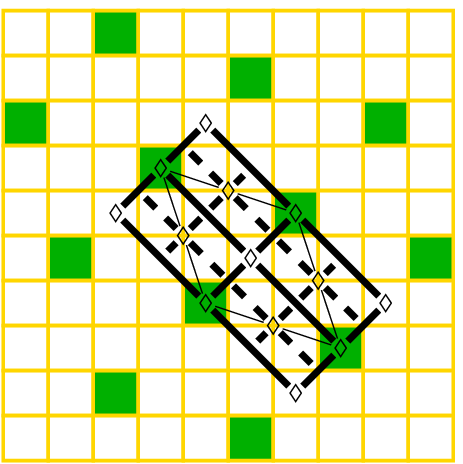} &\epsffile{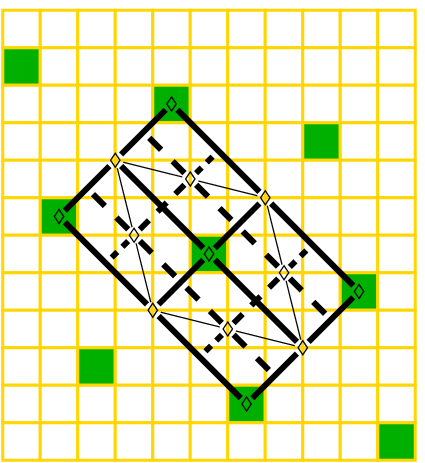}\\
\mbox{(a)} &\mbox{(b)}
\end{array}
\]
\caption{Rhombic satins. a. The (8, 3) satin. \hskip 10 pt b. The (15, 11) satin.}
\label{fig:0a;0b}
\end{figure}

The simplest examples are twills in which each line (vertical or horizontal) of the pattern is formed from the line beside it by an \emph{offset} of one cell.
There are two examples in Fig.~\ref{fig:8d;8e}. 
It is obvious that parallel lines are the same, and it is not hard to see that warps reflect to wefts.
In those diagrams, as in Fig.~\ref{fig:0a;0b}, the broad black lines indicate mirrors of that symmetry in perpendicular directions.
Something to note about oblique reflection is that it reverses which strands are vertical and which are horizontal; so the design as simply reflected has its colours reversed by the colouring convention (vertical dark, horizontal pale).
In order to restore the colours to the original, it is necessary to reverse the colours again with $\tau$.
The broad dashed lines in Figg.~\ref{fig:8d;8e} and \ref{fig:0a;0b} (four in each; four are very short) represent glide-reflections that reverse warp and weft like mirrors and so must be combined with $\tau$; so each dark diagonal band of Figg.~\ref{fig:8d;8e} and \ref{fig:0a;0b} is glide-reflected to the next diagonal band.
Obviously reflections and glide-reflections with horizontal or vertical axes do not reverse warp and weft, but we have none of them in the fabrics we are considering here.

What are in a way the second-simplest examples are satins, where there is very little dark in the pattern (as in the twill examples of Fig.~\ref{fig:8d;8e}, one cell per one-dimensional period, called \emph{order} \cite{Thomas2009}) but where the offset is greater than one moving from strand to adjacent strand.
The examples of Fig.~\ref{fig:0a;0b} are rhombic satins because the dark cells are at the corners of rhombs that are not square.
The first example (Fig.~\ref{fig:0a;0b}a) has offset to the left of three from row to next row up and order 8, and the second (Fig.~\ref{fig:0a;0b}b) has offset four to the left and order 15.
To make the fabric isonemal, the offset has to be relatively prime to the order \cite{GS1980}.

\begin{figure}
\[
\begin{array}{cc}
\epsffile{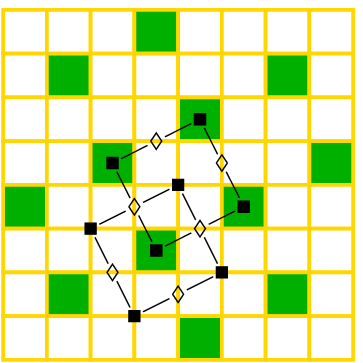} &\epsffile{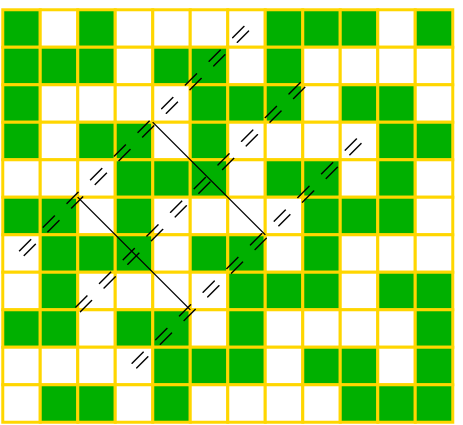}\\
\mbox{(a)} &\mbox{(b)}
\end{array}
\]
\caption{a. The square (5, 3) satin with alternative $G_1$ lattice units outlined. b. Roth's example of species $1_m$, a twillin.}
\label{fig:two;10a}
\end{figure}

These examples have perpendicular mirrors and axes of glide-reflection and so inevitably centres of half-turns where they intersect, represented by diamonds that are hollow because they do not need $\tau$ to be symmetries.
There are a number of ways in which the four examples so far are anomalous.
Their \emph{lattice units}, that is fundamental blocks whose vertices are images of one another under symmetry translations \cite{Schatt1978}, are all rhombs (outlined in the figures), whereas most species of isonemal fabrics have rectangular (Fig.~\ref{fig:two;10a}b) or square (Fig.~\ref{fig:two;10a}a) lattice units.
And so far nothing has had quarter-turn symmetries, and even some satins have that feature, the so-called square satins where the dark cells are the corners of squares.
For example, the satin of Fig.~\ref{fig:two;10a}a has offset 2 to the right ascending and order 5.
The diagram illustrates two different square lattice units, one centred on a dark cell, the other with dark-cell corners.
It also illustrates centres of quarter-turns that need $\tau$ to make them symmetries, represented by filled square boxes \cite{Thomas2010b}.
In this diagram, as in Figg.~\ref{fig:8d;8e} and \ref{fig:0a;0b}, thin lines outline lattice units.

Now that you have seen the attractive pictures to which I have referred and perhaps some others, I define some technical terms.
The fabric designs shown so far have as a symmetry translating a row to the adjacent row and offsetting it by one (a twill) or more cells (a satin or, with more dark cells as in Fig.~\ref{fig:two;10a}b, a twillin).
Prefabrics with those symmetries are said to be of genus I \cite{GS1985}.
Genus II has translation to each adjacent row as for I but combined with  $\tau$.
These genera can overlap, as can the remaining three.
Genus III has rotation to each adjacent row by a half-turn without $\tau$.
Genus IV has rotation to each adjacent row by a half-turn with $\tau$.
Genus V has rotation to the row adjacent in one direction by a half-turn without $\tau$ and to the row adjacent in the opposite direction by a half-turn \emph{with} $\tau$.
You may notice that the overlap of these is illustrated by plain weave, that with diagram a chess board, which lies in all five genera.

\begin{figure}
\[
\begin{array}{cc}
\epsffile{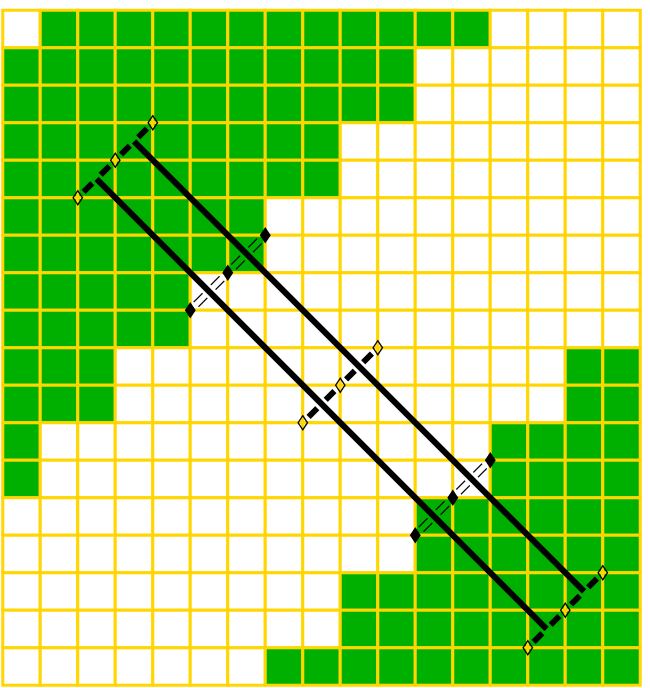} &\epsffile{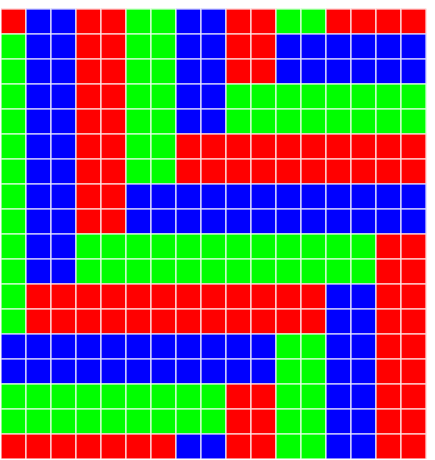}\\
\mbox{(a)} &\mbox{(b)}
\end{array}
\]
\caption{a. Design of a species-$22$ fabric of order 24 with symmetry group partly displayed.\hskip 5 pt b.~Thick striping with three colours, obverse (and reverse).}
\label{fig:6a;6b}
\end{figure}

\begin{figure}
\[
\begin{array}{cc}
\epsffile{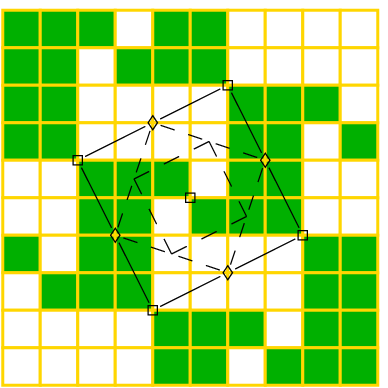} &\epsffile{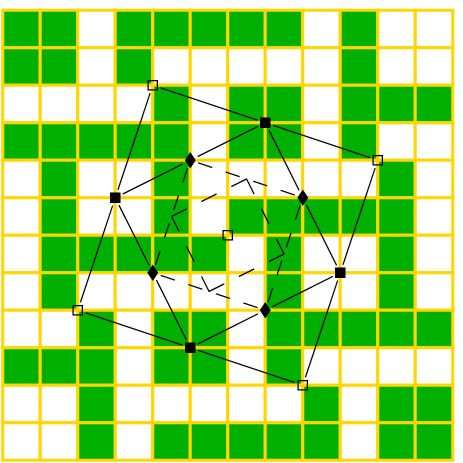}\\
\mbox{(a)} &\mbox{(b)}
\end{array}
\]
\caption{a. Roth's example of species $33_3$, 10-55-2 with $G_1$ outlined and lattice units of levels 1 and 2 dashed (Figure 11b of \cite{Thomas2010b}). b. Roth's species-38 example 20-19437 (Fig.~8 of \cite{Roth1993}) with one lattice unit of $G_1$ and of (larger) $H_1$ outlined and lower-level lattice units dashed (Fig.~12a of \cite{Thomas2010b}).}
\label{fig:81a;85a}
\end{figure}

As the matter of genera illustrates, the symmetry of weaving patterns can be rather complex.
So much so that in the first few years of work on the subject no one took the trouble to sort them out entirely.
This was finally done by Richard Roth \cite{Roth1993}, followed up by a paper \cite{Roth1995} on the two-colouring of the strands not necessarily according to the convention.
Roth found that the \emph{types} of the symmetry groups could be represented well by the pair consisting of the crystallographic type of the group $G_1$ of all the symmetry operations with or without reversal of dark and pale and the crystallographic type of the group $H_1$ of the symmetry operations without reversal of dark and pale.
These pairs represent the types of groups, the latter group $H_1$ being the \emph{side-preserving} subgroup of the full symmetry group $G_1$.
(In Fig.~\ref{fig:two;10a}, the lattice units of $H_1$ and of $G_1$ coincide.)
He writes the pair with the types separated by a stroke as though it represented a factor group.
When the side-preserving subgroup is the whole group, he writes $p4/-$ so that $p4/p4$ means that the subgroup is of the same type $p4$ as the group but is a proper subset.
I'm not going to list them all but, for fabrics of order more than four, there are 39 allowable combinations and so 39 Roth types of symmetry groups of isonemal prefabrics (falling apart doesn't matter to symmetry).
The 39 types can be divided into three broad categories.
1--10 have only parallel axes of symmetry (Fig.~\ref{fig:two;10a}b) \cite{Thomas2009}.
11--32 have perpendicular axes of symmetry (Figg.~\ref{fig:8d;8e}, \ref{fig:0a;0b}, \ref{fig:6a;6b}) \cite{Thomas2010a}.
33--39 have no axes of symmetry but only quarter-turns and half-turns (Figg.~\ref{fig:two;10a}a, \ref{fig:81a;85a}, \ref{fig:7b;85b}b) \cite{Thomas2010b}.
The reason for excluding prefabrics of order 2, 3, and 4 is that, while there are only a few of them, they are mostly exceptions that don't fit into the general schemes.
For instance, some of them have vertical and horizontal axes of symmetry.
Accordingly their symmetry groups are called \emph{exceptional}, and so I'll be discussing what is not exceptional.

I have illustrated most of the possible symmetry operations and their symbols in diagrams. 
It would be a good idea to finish, since the interest in this material lies mainly in how attractive fabric designs often are.
The genus-I twillin 12-183-1%
\footnote{12-183-1 is the catalogue number of this fabric in the second catalogue of isonemal fabrics \cite{GS1986}.
The first catalogue \cite{GS1985} listed non-twill designs up to order 8 plus 13.
Prefabrics that fall apart were catalogued up to order 16 in \cite{HT1991}.}
in Fig.~\ref{fig:two;10a}b illustrates axes of side-preserving glide-reflection represented by hollow dashed lines.
The reversal of dark and pale in the glide-reflection is allowed to stand.
Not illustrated is the other of the two positions that glide-reflection axes can have.
In Fig.~\ref{fig:two;10a}b they pass through the centres and corners of cells, the position that a mirror always has, but they can be differently situated in the only other possible position with respect to cells---through mid-sides, half way between the other positions.
Their position in Fig.~\ref{fig:two;10a}b is called \emph{mirror position} \cite{Thomas2009}.
The  genus-II twillin of Fig.~\ref{fig:6a;6b} illustrates the point that side-reversing and side-preserving glide-reflections can be combined, and furthermore can be combined with mirrors.
It also illustrates the possibility of half-turns that are not side-preserving.
Their centres are represented by filled diamonds.
The term \emph{reverse} indicates the side of the prefabric opposite to the arbitrarily chosen \emph{obverse} side but displayed as though a mirror were held up behind the fabric. This device, which is not original but used in Figg.~185 and 186 of \cite{Emery}, allows one to see how cells are arranged on both sides similarly oriented.
You have seen in Fig.~\ref{fig:two;10a}a quarter-turns that required $\tau$, but quarter-turns do not always require that.
Roth's example 10-55-2 in Fig.~\ref{fig:81a;85a}a, of genus III, illustrates quarter-turns without $\tau$ so that the reversal of dark and pale attendant on a quarter-turn is allowed to stand.
Roth's example 20-19437 in Fig.~\ref{fig:81a;85a}b, of genus IV and V, illustrates that the two sorts of quarter-turns can both appear in a symmetry group but that, when they do, the half-turns at mid-sides of the square lattice units are of the side-reversing type represented by the filled diamonds. 
Dashed lines in Fig.~\ref{fig:81a;85a} will be explained next.

The last two examples illustrate the only other feature of the diagrams that I need to mention.
Square lattice units are highly constrained in their sizes.
At what I call the first level, their sloping sides are the hypotenuses of right triangles with horizontal and vertical sides odd and even and relatively prime.
The smaller dashed squares in Fig.~\ref{fig:81a;85a} are such level-one lattice units with triangle sides 1 and 2, but the fabrics do not have level-one lattice units.
Level-two lattice units are twice that size, also dashed in these diagrams; I draw them origami-style so that the corners of the smaller are mid-sides of the next larger.
The $G_1$ lattice units of these designs are at level three, twice the size of the level-two dashed squares inside them.
Level four is illustrated in Fig.~\ref{fig:81a;85a}b, with the lattice unit of the side-preserving subgroup, the outermost square.
Fig.~\ref{fig:7b;85b}b is a design with a larger lattice unit at level four for its symmetry group; it is its own side-preserving subgroup.
In the side-preserving subgroup of 20-19437 (Fig.~\ref{fig:81a;85a}b), because the filled diamonds and filled boxes of the symmetry group, involving $\tau$, cannot appear in the side-preserving subgroup, the diamonds simply disappear and the filled boxes at mid-sides should appear as hollow diamonds, since the composition of two side-reversing quarter-turns is a side-preserving half-turn.
What's left is the side-preserving quarter-turns and the half-turns where the other quarter-turns used to be.
The lattice unit is bigger than that for $G_1$.
Finally in this discussion of lattice units, the centre of a square lattice unit can be either at corners of cells as in Figures 3a (upper) and \ref{fig:81a;85a} or at the centre of a cell, as in the lower example in Fig.~\ref{fig:two;10a}a.
In the former case, lattice-unit levels can go only as far as four, and in the latter case only as far as two \cite{Thomas2010b}.

Once one knows what the groups are, one knows the sizes of the lattice units, and it's a matter of filling in the blanks to determine all of the isonemal fabrics with lattice units of any particular size.
I found in \cite{Thomas2009}, \cite{Thomas2010a}, and \cite{Thomas2010b} that I had to refine Roth's taxonomy of groups somewhat to produce a taxonomy of fabrics.
So where he has types of group numbered from 1 to 39 with many inclusions---groups that are subgroups of other groups, I have species and subspecies of prefabrics that do not overlap numbered $1_m$, $1_o$, $1_e$, $2_m$, $2_o$, $2_e$, 3, $4_o$, $4_e$, $5_o$, $5_e$, 6, \dots $35_3$, $35_4$, $36_1$, $36_2$, $36_s$, 37, 38, 39, where the group of a prefabric of species $n_x$ is of type $n$ and $x$ indicates a disjoint subdivision, if any.
For example $1_m$ indicates that glide-reflection axes are in mirror position as in Fig.~\ref{fig:two;10a}b, whereas $1_o$ and $1_e$, while different from each other, both have axes not in mirror position.
Numerical subscripts indicate levels of square lattice units, and $s$ indicates a subset of the level-two designs containing, \emph{inter alia}, the square satins of even order.
The square satins of odd order (Fig.~\ref{fig:two;10a}a) fall in subspecies $36_1$.

Before the weaving of cubes was considered in \cite{Thomas2010b} and their 2-colouring in \cite{Thomas2012}, the analogous weaving of tori should have been mentioned. A standard weaving design is a 2-colour tessellation of the plane, and such a 2-colour tessellation can be interpreted as a weaving design for a flat torus in several ways, although not arbitrarily. The {\em locus classicus} of flat torus tessellation is \cite{Coxeter} building on \cite{Brahana}. It makes no sense to identify as opposite sides of the map to the torus what is not identical in the planar design. For a design of order $n$, the simplest thing to do is to use an $n$-by-$n$ square of cells as the 2-dimensional period parallelogram and map it to the torus. For example, a quarter of Fig.~\ref{fig:8d;8e}a. There are then $n$ strands in the vertical and the horizontal directions; isonemality makes sense and is preserved. If an $n$-by-$2n$ rectangle of cells were used, opposite sides could be identified, but horizontal and vertical strands could not be interchanged by an isometry. On the other hand, an $mn$-by-$mn$ square would do as well as making $m=1$. An example with $m=2$ is the whole of Fig.~\ref{fig:8d;8e}a. Perhaps more interesting is to use an $H_1$ lattice unit as period parallelogram or, when it is a rhomb, the rectangle containing it as within Fig.~\ref{fig:8d;8e}a. The number of strands in each direction, now oblique, will be reduced by the identification, in the case of Fig.~\ref{fig:8d;8e}a to a single strand. The reduction in the number of strands occurred also in the weaving of cubes in \cite{Thomas2010b}, where a design of order 10 required only 6 strands to weave the cube using an oblique $G_1$ lattice unit for each face and the order-20 design of Fig.~16b required only 8 strands. The latter design (topological, not visual) is used in Figg.~14 and 15 of \cite{Thomas2012} to 2-colour a woven cube. This is not the place to discuss what the definition of an isonemal weaving of a torus or Klein bottle should be. 

\section{Unlimited colours}
\label{sect:Unlimited}

\noindent The visual appearance of a coloured fabric from the obverse side I refer to as a \emph{pattern}.
Such a pattern, when it consists of an array of dark and pale congruent cells tessellating the plane, is given a topological meaning, the \emph{design} of a fabric.
As we shall be considering patterns different from the design of a fabric, the distinction is important; a design is the pattern of a fabric that is normally coloured.
So a two-colour pattern can be interpreted as a design or not.
Normal colouring has a consequence that other colourings of the strands need not have.
Because at every point not on the boundary of a strand the two strands preferentially ranked are of two different colours, a design's colour complement (switching dark and pale) is the appearance of the fabric from behind as though viewed from the obverse side in a mirror set up behind $E$, the plane of the fabric, its reverse.
When a fabric is coloured normally the reverse pattern is the colour complement of the obverse pattern.

\begin{figure}
\[
\begin{array}{ccc}
\epsffile{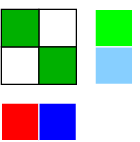} &\raisebox{16 pt}{\epsffile{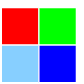}} &\raisebox{16 pt}{\epsffile{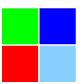}}\\%
\mbox{(a)} &\mbox{(b)} &\mbox{(c)}
\end{array}
\]
\caption{a. Design of plain weave and one choice of colours for a thin striping. \hskip 10 pt b. Obverse view of colouring (a). \hskip 10 pt c. Reverse view of colouring (a).}
\label{fig:1a:1b;1c}
\end{figure}

\begin{figure}[t]
\[
\begin{array}{ccc}
\epsffile{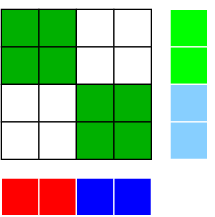} &\raisebox{16 pt}{\epsffile{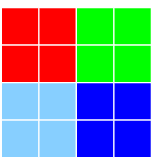}} &\raisebox{16 pt}{\epsffile{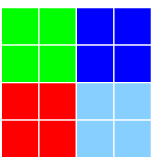}}\\
\mbox{(a)} &\mbox{(b)} &\mbox{(c)}\\
\\
\epsffile{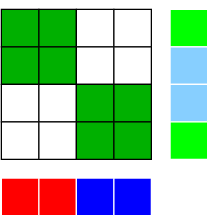} &\raisebox{16 pt}{\epsffile{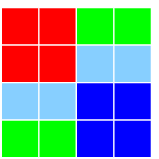}} &\raisebox{16 pt}{\epsffile{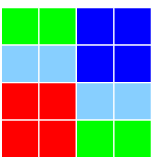}}\\
\mbox{(d)} &\mbox{(e)} &\mbox{(f)}\\
\\
\epsffile{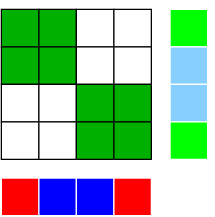} &\raisebox{16 pt}{\epsffile{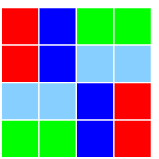}} &\raisebox{16 pt}{\epsffile{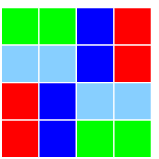}}\\
\mbox{(g)} &\mbox{(h)} &\mbox{(i)}
\end{array}
\]
\caption{Plain weave doubled.\hskip 10 pt a. Colouring by thick striping to preserve the blocks of four. \hskip 10 pt b. Obverse view. \hskip 10 pt c. Reverse view.\hskip 10 pt d. A second colouring by thick striping. \hskip 10 pt e. Obverse view. \hskip 10 pt f. Reverse view.\hskip 10 pt g. A third colouring by thick striping. \hskip 10 pt h. Obverse view. \hskip 10 pt i. Reverse view.}
\label{fig:2a;2b;2c;3a;3b;3c;4a;4b;4c}
\end{figure}

To this point colours (dark and pale) have been used for their topological meaning.
They are now going to be used more generally.
We adopt a device due to Roth to indicate, when not obvious, how strands are coloured, thinly in Fig.~\ref{fig:1a:1b;1c} and thickly in Fig.~\ref{fig:2a;2b;2c;3a;3b;3c;4a;4b;4c}.

A symmetry operation of a fabric with coloured strands is called a \emph{colour symmetry} \cite{Roth1995} `if it permutes the colors consistently'.
All the strands of one colour must be mapped either to strands of that colour or to strands of some other colour, and correspondingly the strands of each colour.
If all of the weave symmetries of a fabric with coloured strands are colour symmetries, then the choice of the strand colours is said to be \emph{perfect} or \emph{symmetric}, where I shall use exclusively the former term.
The interaction of design and pattern is important because the weave symmetry group is that of the fabric represented by the two-colour pattern that is its design, the `group of the design' for short, but the permutation (identity or other) of colours occurs only in the pattern that represents the colouring of the strands, identical to the design only with normal colouring of the strands.
There is potential for confusion, especially when colouring with two colours.

Perfect colouring with two colours is one subject of Roth's second weaving paper \cite{Roth1995}.
The so-called \emph{normal colouring} that produces designs of woven fabrics is a perfect colouring.
When one seeks other perfect two-colourings \cite{HT1991}, one finds that one can colour the strands only alternately in both directions in the two colours or one can colour the strands in adjacent pairs in both directions in the two colours.
These are called \emph{thin} and \emph{thick} striping of the warp and weft.
In his paper Roth determines which fabrics---actually prefabrics---can be perfectly coloured by striping warp and weft with two colours but not in terms of his taxonomy.
It has been done in those terms and some further results deduced in \cite{Thomas2011} and \cite{Thomas2012}.
This paper is intended to extend that determination to colouring with more than two colours.

Thin stripes overlap in individual cells---all cells---of the fabric, but thick stripes overlap in square blocks of four cells, which will be referred to as \emph{blocks}.
The natural unit of measurement in the fabric is the side of a cell, in which the diagonal of the square cells have the length $\sqrt 2$, which will be abbreviated $\delta$.
We begin by exploring thick striping without limiting the number of colours.

\begin{Lem} 
The assignment of all different colours to a thick striping of a non-exceptional isonemal fabric is a perfect colouring of the fabric provided that the thick stripes are preserved by the symmetries of the fabric, i.e., that

\noindent 1. components of translations, both horizontal $x$ and vertical $y$ are even in cell widths;

\noindent 2. centres of half-turns are confined to corners of cells;

\noindent 3. centres of quarter-turns are confined to corners and centres of blocks;

\noindent 4. diagonal glide-reflections have either

\noindent a. axes through block centres and glides even in $\delta$ or

\noindent b. axes through the centres of cells but not through block centres and glides odd in $\delta$;

\noindent 5. diagonal axes of reflection pass through block centres.
\label{lem:1}
\end{Lem}

\noindent \emph{Proof.} The restrictions on the symmetries of the fabric are just those that preserve the configuration of blocks.
Since all of the stripes have different colours, any symmetry operation that preserves the blocks will be a colour symmetry.
Colour simply does not enter into consideration beyond defining the blocks.\hfill\whbox

The restrictions of Lemma \ref{lem:1} are an easing of the restrictions for thick striping with two colours because with duplication of warp and weft colours there are \emph{redundant} blocks where stripes of the same colour cross.
Perfect two-colouring requires not only the preservation of the configuration of the blocks but the preservation of the configuration of redundant blocks, which form a checkerboard with their complement, the \emph{irredundant} blocks.
(Specifically, the additional restrictions for 2-colourings rule out operations that exactly interchange the redundant and the irredundant blocks:

\noindent 1. The even $x$ and $y$ must not be congruent to 0 and 2 or 2 and 0 (mod 4).

\noindent 2. Centres of half-turns must not be at mid-sides of blocks.

\noindent 3. Centres of quarter-turns must not be at corners of blocks.

\noindent 4. Diagonal glide-reflection axes must not miss block centres (\emph{a} above allowed, \emph{b} not).

\noindent There is no further restriction on mirrors.)

Even for arbitrarily many colours, the translation constraint immediately eliminates strand-to-adjacent-strand translations (genera I and II), taking with them species 1--10, 12, 14, 16, 18, 20, 23, 24, 26, 28, 31, 32, 34, 36, and 39.
Call this the \emph{twillin ban}.
On the other hand, we know from \cite{Thomas2012} that fabrics of species 17, 19, 21, 25, 27, 29, 33, 35, and 37 can be thickly striped with only two colours, a more tightly constrained task.
We need to examine only species 11, 13, 15, 22, 30, and 38 to see which can have strands striped thickly if the number of colours is unrestricted.

In order to ignore differences that do not make a difference, it is necessary to introduce a new group that ignores reflections $\tau$ in the plane of the fabric. 
Every element of a symmetry group is of the form $(t, r)$, where $t$ is a transformation in the plane and $r$ is either $\tau$ or the corresponding identity $e$.
For the present purpose we need the projection $G_1 \rightarrow G_2$ onto the group $G_2$ consisting of elements of the form $(t, e)$ for all $(t, r)\in G_1$ whether $r=\tau$ or $r=e$.
$H_1$ is quite different for some fabrics, being the subgroup of elements of $G_1$ that are of that form; when $H_1 = G_1$, $G_2$ is also $G_1$.
When $H_1$ differs from $G_1$, it will also differ from $G_2$ since each element of $G_1$ not in $H_1$ will appear in $G_2$ stripped of its $\tau$.
In the example of Fig.~\ref{fig:8d;8e}a, the reflections and glide-reflections, including $\tau$ as they do, are omitted from $H_1$ and so appear differently in $G_2$, where they are side-reversing reflections and glide-reflections rather than side-preserving.

\begin{figure}
\[
\begin{array}{ccc}
\epsffile{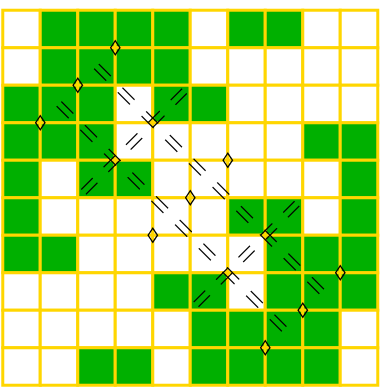} &\epsffile{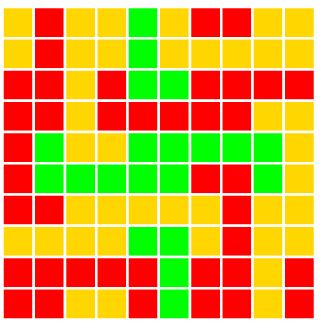} &\epsffile{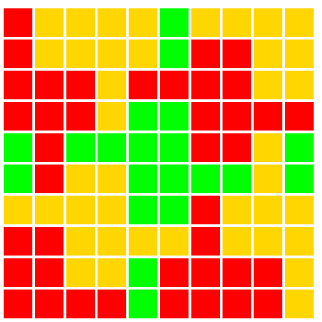}\\
\mbox{(a)} &\mbox{(b)} &\mbox{(c)}
\end{array}
\]
\caption{a. 12-111-2 of subspecies $11_o$.\hskip 10 pt b. Thickly striped with 3 colours, obverse. \hskip 10 pt c. Reverse.}
\label{fig:5c;5a;5b}
\end{figure}

Species 11, 13, and 15 with perpendicular glide-reflection axes differ only in the use of $\tau$, the reflection in the plane $E$, not within the side-preserving $G_2$; so they can be discussed together.
Subspecies $11_o$, $13_o$, and $15_o$ have all of their glides odd and present no difficulty in positioning glide-reflection axes with respect to the blocks.
An example of a perfect colouring of Roth's example 12-111-2 of subspecies $11_o$ with three colours (a more constrained task than with all stripes different colours) is shown in Fig.~\ref{fig:5c;5a;5b}.

Species 22 has odd glides perpendicular to its mirrors.
The glide-reflection axes must, by 4b of Lemma \ref{lem:1}, not pass through block centres, and so the centres of half-turns, lying on the glide-reflection axes, must fall at mid-sides of blocks.
This species' constraints present no difficulty in fitting axes to blocks.
Again a 3-colour example illustrates the point in Fig.~\ref{fig:6a;6b}b with a design made up for the purpose (Fig.~\ref{fig:6a;6b}a).
Because of the mirror symmetry, which on account of involving $\tau$ appears only in the relation between obverse pattern and reverse pattern, one needs to see both sides to check the mirror symmetry.
On the other hand, its presence makes the reverse have the same pattern as the obverse and so be of no distinct visual interest.
Only the diagonal lines of redundant blocks are identical on the two sides of the fabric; the reverse has horizontal stripes approximately where the obverse has vertical stripes, and vice versa---only approximately because, as is easily seen, the stripes do not end exactly at the redundant blocks.
Such reverse patterns will not be illustrated.

\begin{figure}
\[
\begin{array}{cc}
\epsffile{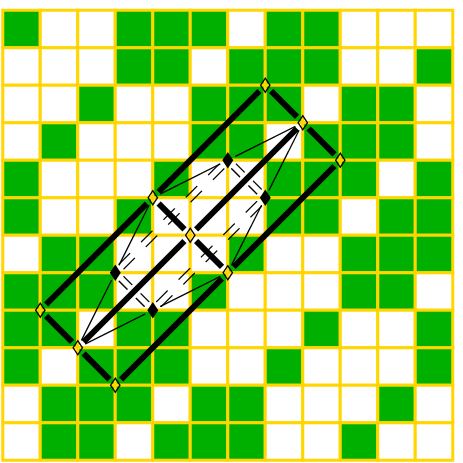} &\epsffile{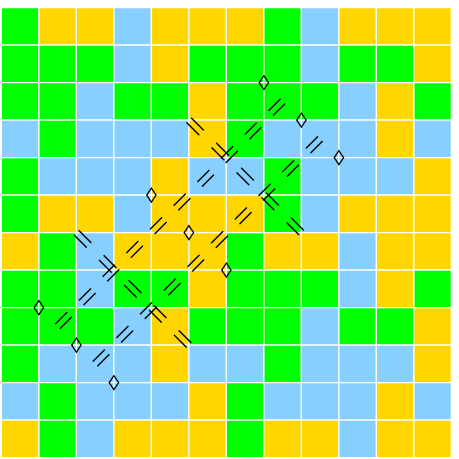}\\
\mbox{(a)} &\mbox{(b)}
\end{array}
\]
\caption{a. Design of Roth's example 12-315-4 of species 30. b. Pattern of thick striping with three colours and no mirror symmetry.}
\label{fig:13a;7a}
\end{figure}

Species 30, with the spacing of subspecies $27_o$ and not perfectly 2-colour\-able with thick stripes, has odd glides like $11_o$, $13_o$, and $15_o$ and presents no difficulty to the use of more colours than two.
For example, Roth's example 12-315-4 (Fig.~\ref{fig:13a;7a}a) can be thickly striped with three colours (Fig.~\ref{fig:13a;7a}b).
This is the first example with interesting motifs having more symmetry than the pattern as a whole; another perfect colouring of this fabric is a herring-bone pattern like Fig.~\ref{fig:6a;6b}b but of course with shorter stripes.
The symmetry of the pattern of Fig.~\ref{fig:13a;7a}b is notable not only because of the quarter-turn-symmetric houndsteeth linked/separated by dominos. 
Each colour as a whole has rotational symmetry---not shared with the other two colours---with centres of quarter-turns in the centres of that colour's houndsteeth, where the half-turn centres marked in Fig.~\ref{fig:13a;7a}b apply to the whole pattern (indeed to the design).
As in species 22, the mirror symmetry is revealed only in the relation between obverse and reverse which, because of those mirrors and side-reversing half turns, looks just like the obverse.

\begin{figure}[t]
\[
\begin{array}{cc}
\epsffile{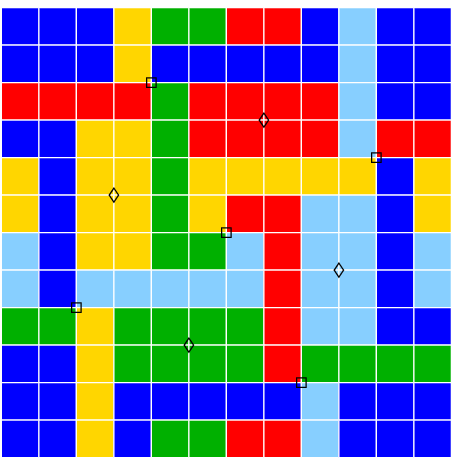} &\epsffile{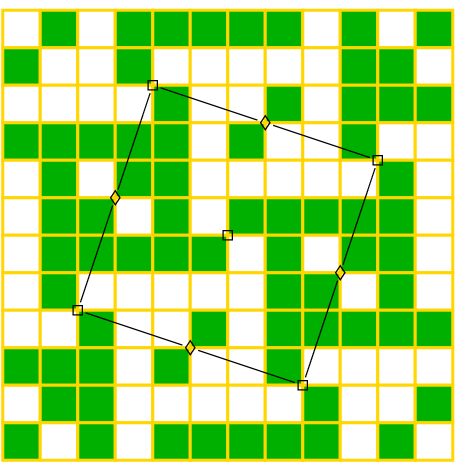}\\
\mbox{(a)} &\mbox{(b)}
\end{array}
\]
\caption{a. Roth's 20-19437 of species 38 (Fig.~\ref{fig:81a;85a}b) thickly striped with five colours. b. Order-20 example (Fig.~13 of \cite{Thomas2010b}) of species $33_4$ with $G_1$ of level 4 outlined.}
\label{fig:7b;85b}
\end{figure}

Species 38 has no symmetry axes but has level-3 lattice units and both side-preserving and side-reversing quarter-turns (for example, Fig.~\ref{fig:81a;85a}b).
If blocks are centred on one sort of quarter-turn centres, then the level guarantees that the other sort of quarter-turn centres will fall at the corners of blocks.
The level also guarantees that translation components will be even as required.
While examples with three and four colours are impossible, as will be shown at the end of the next section (the quarter-turn ban), Roth's example 20-19437 (Fig.~\ref{fig:81a;85a}b) can be perfectly coloured with five colours (Fig.~\ref{fig:7b;85b}a).
The visible (and marked) symmetry of the figure is that of the design's $H_1$ at level 4, the side-reversing quarter-turns of $G_1$ appearing only squared as half-turns, otherwise disappearing from a one-sided point of view.
Because of the side-reversing quarter-turns and half-turns, the reverse is the same as the obverse and so not illustrated.
In contrast with Fig.~\ref{fig:13a;7a}b, the motifs of Fig.~\ref{fig:7b;85b}a have less symmetry than the pattern, the symmetry of which arises from how the motifs are arranged.
 
The result of the above exploration is that some fabrics in all species not eliminated by the twillin ban (conflict of their genus with the translation constraint) are perfectly colourable by thick stripes of all different colours, indeed even by as few as five colours, which is what has actually been illustrated.

\begin{Lem} 
The assignment of different colours to all strands of an isonemal fabric is a perfect colouring of the fabric by thin striping.
\label{lem:2}
\end{Lem}

\noindent \emph{Proof.} All of the limitations to preserve thick stripes are automatic in their thin-stripe analogues.
All symmetries preserve cells, hence strands, and that is all that is needed when all the strands have different colours.\hfill\whbox

\begin{Thm} 
Some fabrics in all non-exceptional species can be perfectly colour\-ed by the assignment of different colours to all single strands or, with the exception of fabrics in genus I or II, to all adjacent pairs of strands.
\label{thm:1}
\end{Thm}

\noindent \emph{Proof.} Lemma \ref{lem:2} for thin striping and the discussion between Lemma \ref{lem:1} and Lemma \ref{lem:2} for thick striping show this.\hfill\whbox

\section{A finite number of colours}
\label{sect:Finite}

\begin{Lem} 
The assignment of a finite number of colours to weft and warp stripes of an isonemal fabric can be a perfect colouring only if the same number of colours are assigned to warps as to wefts.
\label{lem:3}
\end{Lem}

\noindent \emph{Proof.} Since the fabric is isonemal, there is a symmetry operation taking each warp to each weft.
Under such an operation, if it is a colour symmetry, each warp colour is mapped to a weft colour.
So the number of colours of warps is the number of colours of wefts.\hfill\whbox

\begin{Lem} 
The assignment of a finite number of colours to weft and warp stripes of an isonemal fabric can be a perfect colouring only if the vertical and horizontal sequences of the colours of the stripes are periodic with the same period.
\label{lem:4}
\end{Lem}

\noindent \emph{Proof.} Parallel stripes of the same colour must be the same distance from their same-coloured parallel neighbours because each one can be mapped to each one by a colour symmetry that acts as the identity on their colour.
But since the strands of each colour can also be mapped to the parallel strands of any other colour also by a colour symmetry, that inter-stripe distance is constant over the colours horizontally and vertically.
And since the warps can be mapped to the wefts by a colour symmetry, the interstripe distance or period is common to warps and wefts.\hfill\whbox

What makes achieving perfection in colouring by striping difficult rather than impossible is redundancy --- where strands cross strands of the same colour.
Using only a finite number of colours, as section \ref{sect:Unlimited} suggested, is not a difficulty.

\begin{Thm} 
If no colour is common to warps and wefts in the assignment of a finite number of colours to the striping of an isonemal fabric, then the necessary conditions of Lemmas \ref{lem:3} and \ref{lem:4} allow perfect colouring by thin striping and the necessary conditions of Lemmas \ref{lem:1}, \ref{lem:3}, and \ref{lem:4} allow perfect colouring by thick striping.
\label{thm:2}
\end{Thm}

\noindent \emph{Proof.} Obvious.\hfill\whbox

\begin{Lem} 
The assignment of a finite number of colours to stripes of an isonemal fabric can be perfect only if the colours of the warps are the colours of the wefts or no colour is shared.
\label{lem:5}
\end{Lem}

\noindent \emph{Proof.} If any colour $C$ is common to warps and wefts, then let colour $D$ be any other weft colour. Since the colouring is perfect, there is a colour symmetry mapping the wefts of colour $C$ to wefts of colour $D$, but also mapping warps of colour $C$ to warps of colour $D$.
But then there must be warps of colour $D$.
Each weft colour is also the colour of some warps.
And, since there are the same number of colours of warps and wefts, the colours are the same.\hfill\whbox

While the warp colours, if not all different, must be the same as the weft colours, examples show that the order of the warp colours need not be the same as the order of the weft colours.
One can see in Fig.~\ref{fig:7b;85b}a that there are red and pale blue wefts adjacent to yellow wefts but not red or pale blue warps adjacent to yellow warps.
When the $c$ colours, say, of the warps are the colours of the wefts, then in a $c\times c$ square of cells/blocks, for thin/thick striping respectively, each warp colour crosses each weft colour in exactly one cell/block.
So there are exactly $c$ redundant cells/blocks, one in each row and column of the $c\times c$ square of cells/blocks.
In order for the colouring to be a perfect colouring, the symmetry group of the fabric must leave the redundant cells/blocks invariant as a whole and each element of the group must permute them, if at all, in a consistent way. 

From now on, warps and wefts will be coloured with the same finite set of colours.
For a specific colouring of the wefts, say, as illustrated in Fig.~\ref{fig:7b;85b}a, the positions of the redundant cells must be chosen in accordance with the consistency requirement so that the symmetry group of the fabric is appropriately related to the symmetry group of the pattern of redundant blocks.
In this case the doubled (5, 3) satin (10-3-2) is a choice of redundant-block distribution that works.
For the design of Fig.~\ref{fig:81a;85a}b and the placement of redundant blocks surrounding the four centres of quarter-turns (extended by the four corner blocks), the downward sequence of colours, dark blue, red, yellow, pale blue, green, (dark blue,) dictates the left-to-right sequence, dark blue, yellow, green, red, pale blue, (dark blue).
This topic cannot be pursued here.

We now state Roth's Colouring Theorem.
\begin{Thm} 

If the colours of the warp stripes and the colours of the weft stripes are the same finite set of colours, then, in addition to the necessary conditions of Lemmas \ref{lem:3} and \ref{lem:4} for thin striping and in addition to the necessary conditions of Lemmas \ref{lem:1}, \ref{lem:3}, and \ref{lem:4} for thick striping, the redundant cells/blocks must collectively be preserved by $G_1$ of an isonemal fabric for its colouring to be perfect.
\label{thm:3}
\end{Thm}

\noindent \emph{Proof.} A symmetry operation is a colour symmetry exactly when it preserves and permutes the colours of the redundant cells/blocks because each brings with it its wefts and its warps, which have the same colour.\hfill\whbox

This theorem was observed to be the case by Roth \cite{Roth1995} for two colours, but the essence is independent of the number of colours.
Colouring the set of redundant cells/blocks specifies a colouring by striping warps and wefts with the same colours.
The requirement that the symmetry group of the fabric not affect the redundant cells can be stated by saying that $G_2$ of the fabric to be striped must be a subgroup of $G_2$ of the pattern of redundant cells.

\begin{Lem} 
The redundant cells/blocks of a perfect colouring of the warps and wefts of an isonemal fabric by thin/thick striping with $c$ colours are the dark cells/blocks of an order-$c$/doubled order-$c$ derived isonemal fabric with one/two dark cells per order length, $c/2c$.
\label{lem:6}
\end{Lem}

\noindent \emph{Proof.} Since a stretch of $c$ cells/blocks along a strand/pair of some colour meets strands/pairs of $c$ different colours, only one of the $c$ cells/blocks is redundant. 
Let the redundant cells/blocks be the dark cells of a derived fabric. 
By Lemma \ref{lem:4}, the colouring is periodic with period $c/2c$. 
The order of the derived fabric is $c/2c$, which must divide the order of the coloured fabric.
The group $G_2$ of the derived fabric has $G_2$ of the coloured fabric as a subgroup because that $G_2$ leaves the dark cells invariant collectively.
Since $G_1$ of the coloured fabric maps every warp and weft to every warp and weft preserving redundant cells/blocks, the derived fabric is isonemal under the action of the coloured fabric's group $G_1$, \emph{a fortiori} under its own symmetry group.\hfill\whbox

For thin striping, the redundant cells may be arranged as $(c-1)/1$ twills, and for $c\geq 4$ as various other isonemal fabrics,
 6-1-1, 8-1-1, 8-1-2, and including in particular the square satins beginning with 5-1-1.
For thick striping, the redundant cells can be arranged as doubled $(c-1)/1$ twills, 6-3-1, 8-3-1, 10-3-1, \dots and various other doubled fabrics beginning with 10-3-2 (doubled 5-1-1), and the doubles of 6-1-1, 8-1-1, 8-1-2, and so on.

It is easy to see that any rigid motion of a line of twilly redundant cells to itself is a coherent permutation of their colours whether the motion merely translates or reverses direction with a cell fixed or not.
Also that the effect on the parallel lines of redundant cells has the same colour consequences.

There is reason to consider mostly unexceptional isomenal fabrics, i.e.,  with order more than four, but in order to consider colourings with three colours, one needs to begin thin striping with the 2/1 twill as a redundancy pattern, there being no other choice.
Because side-reversing symmetries (mirrors and side-reversing glide-reflections, half-turns, and quarter-turns) of a fabric relate its opposite sides and one sees fabrics one side at a time, fabrics with exclusively side-reversing symmetries are less appealing when striped.
Emphasis will be on fabrics with side-preserving symmetries just because their symmetries are visible when the fabric is striped.%
\footnote{While the presence of side-reversing symmetry makes it uninteresting and therefore unnecessary to look at both sides of such coloured fabrics, the absence of all such symmetry makes the two sides of a fabric that lacks it matters of independent interest. E.g., Fig.~\ref{fig:5c;5a;5b}.}

For three and four colours and order greater than four there is no alternative to colouring as one could do with any odd or even number of colours where the order of the colours is the same for warps as for wefts specified by a twill arrangement of redundant cells or blocks.
It is obvious in both cases that no fabric with quarter-turn symmetry can have $G_2$ a subgroup of $G_2$ of the 2/1 or 3/1 twill---or of any other $(c-1)/1$ twill doubled or not, since these twills fall into species with no quarter-turns in their $G_2$s.
(The 2/1 twill is of subspecies $28_o$, the 3/1 twill is of $26_e$, 6-3-1 is of $27_o$, and 8-3-1 is of $25_e$.)
Accordingly, no fabric of species 33--39 can be perfectly coloured with three or four colours---or, with twilly redundancy, any number of colours greater than two.
Call this the \emph{quarter-turn ban} to go with the twillin ban of section \ref{sect:Unlimited}.
\begin{figure}[h]
\epsffile{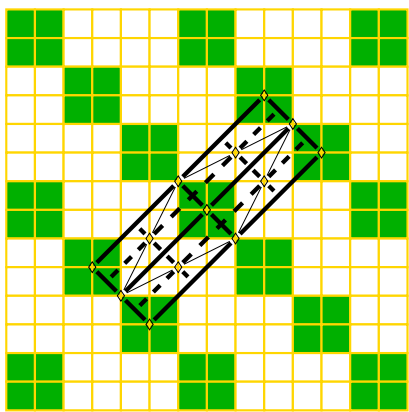}
\caption{Doubled 2/1 twill as the redundant blocks for thick striping with three colours.}
\label{fig:61}
\end{figure}

\section{Three colours}
\label{sect:three colours}

\noindent As an example of how thick striping can be effected to create perfect colouring, the case of three colours will be worked out in this section with twilly redundancy.
The redundant doubled twill for thick striping with three colours 6-3-1 (Fig.~\ref{fig:61}) is of subspecies $27_o$, being of course twill 2/1 doubled.
Species 33--39 having been eliminated by the quarter-turn ban, the twillin ban (section \ref{sect:Unlimited}) reduces the species among which it might be possible to stripe thickly with any finite number of colours for species 1--32 to 11, 13, 15, 17, 19, 21, 22, 25, 27, 29, and 30.
Fabrics for all these species can be perfectly coloured by thickly striping with three colours.

For species 11, 13, and 15, glide-reflection axes that can be an odd multiple ($\geq 1$) of $3\delta$ apart can be placed on the mirrors in the dark blocks (of 6-3-1) if their glides are even ($11_e$, $13_e$, $15_e$) because the blocks along their line can be related by even glides or half-way between them if their glides are odd ($11_o$, $13_o$, $15_o$) because adjacent line's blocks are related by odd glides.
Since the fabric's glides perpendicular to the line of dark blocks are odd in length, their axes can be placed on the glide-reflection axes $\delta$ apart perpendicular to the dark lines.
An example is 12-111-2 of species $11_o$ coloured by thick striping with 3 colours already shown in Fig.~\ref{fig:5c;5a;5b}. 

\begin{figure}
\[
\begin{array}{cc}
\epsffile{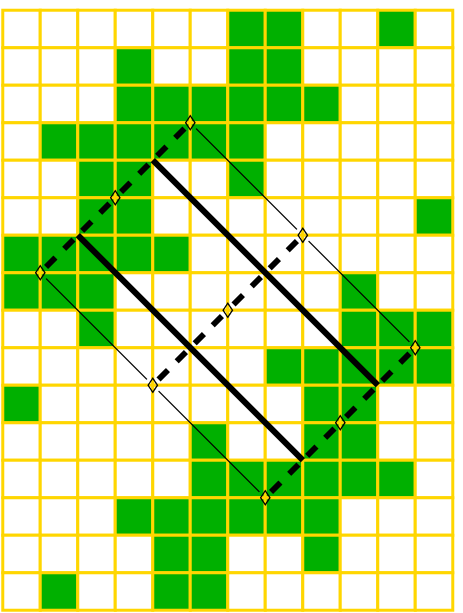} &\epsffile{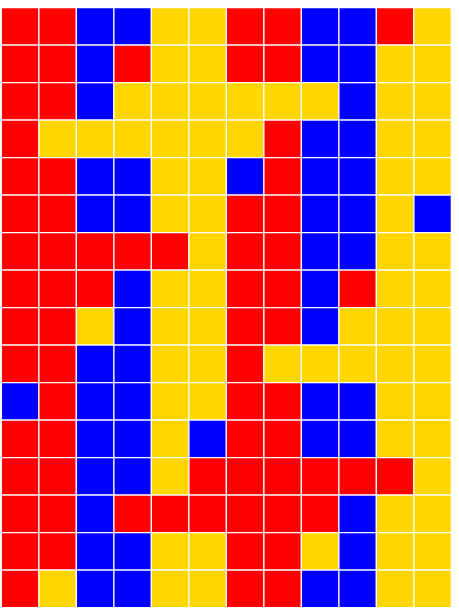}\\
\mbox{(a)} &\mbox{(b)}
\end{array}
\]
\caption{a. Order-24 fabric of species $17_e$. b. Reverse of fabric 3-coloured by thick striping. Obverse, being the reflection, has horizontal stripiness.}
\label{fig:58c;58d}
\end{figure}

\begin{figure}
\[
\begin{array}{cc}
\epsffile{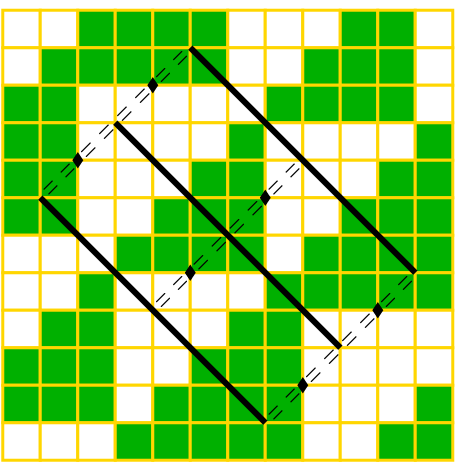} &\epsffile{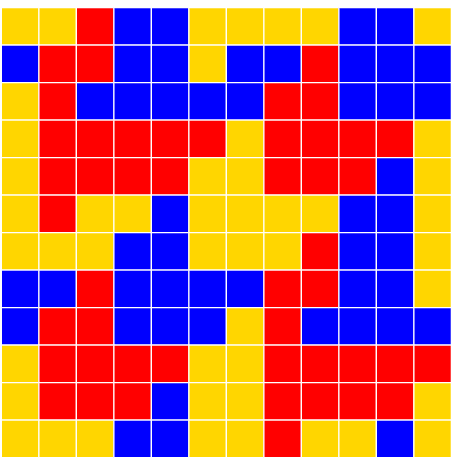}\\
\mbox{(a)} &\mbox{(b)}
\end{array}
\]
\caption{a. Order-24 fabric of species $19_e$. b. Fabric 3-coloured by thick striping. The middle of the figure makes the SW-NE period look like $2\delta$, but the upper left matching yellow and blue motifs make it clear that it is $4\delta$.}
\label{fig:59c;59d}
\end{figure}

\begin{figure}
\[
\begin{array}{cc}
\epsffile{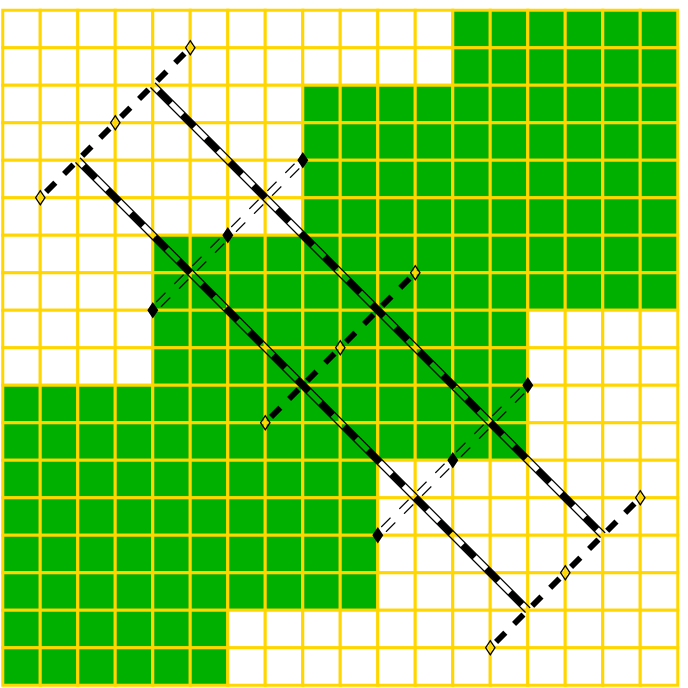} &\epsffile{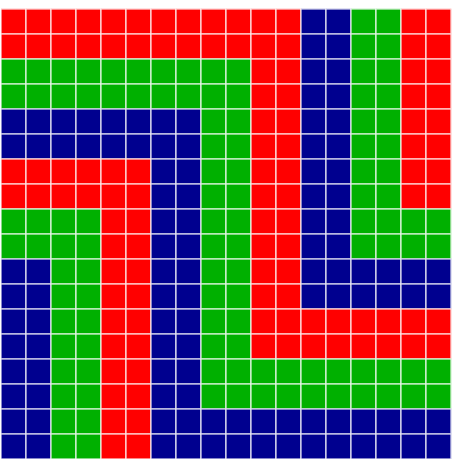}\\
\mbox{(a)} &\mbox{(b)}
\end{array}
\]
\caption{a. Design of a fabric of species 21 and order 96. b. Three-colouring by thick striping with lines of redundant blocks along axes with positive slope.}
\label{fig:18a;62_3}
\end{figure}

For subspecies $17_e$ and $19_e$, and species 21, the distance between the mirrors, which is the glide of the glide-reflections, can be taken to be even.
The axes of the glide-reflections, if an odd multiple ($\geq 1$) of $3\delta$ apart, can therefore be placed on the dark lines' mirrors in the dark blocks.
The mirrors of subspecies $25_e$, if an odd multiple ($\geq 1$) of $3\delta$ apart can, like the axes of glide-reflections with even glides in $17_e$, $19_e$, and 21, go along the mirrors in the dark blocks.
The fabric's other mirrors can then be placed on the mirrors (of 6-3-1) perpendicular to them.
Examples are a fabric of species $17_e$ illustrated in Fig.~\ref{fig:58c;58d}a and coloured in Fig.~\ref{fig:58c;58d}b, a fabric of species $19_e$ illustrated in Fig.~\ref{fig:59c;59d}a and coloured in Fig.~\ref{fig:59c;59d}b, and the fabric of Fig.~\ref{fig:18a;62_3}a of species 21 coloured in Fig.~\ref{fig:18a;62_3}b.

For subspecies $17_o$ and $19_o$ and species 22 the distance between mirrors, which is the glide of the glide-reflections, can be taken to be odd.
The axes of the glide-reflections, if an odd multiple ($\geq 1$) of $3\delta$ apart, can therefore be placed half-way between the dark lines' mirrors along the dark blocks.
The fabric's mirrors then can be placed on the mirrors perpendicular to them.
Examples are Roth's example 12-11-2 of species $17_o$, 12-55-5 of species $19_o$, 
and the fabric of species 22 displayed with its thick 3-colouring in Fig.~\ref{fig:6a;6b}b.

For subspecies $27_e$ and species 29 with the same spacing, glides in both perpendicular directions are even.
Axes that can be an odd multiple ($\geq 1$) of $6\delta$ apart and parallel mirrors can all be placed along the mirrors through the dark blocks.
The fabric's perpendicular axes and mirrors can then be placed along the mirrors perpendicular to them.
\begin{figure}
\noindent
\epsffile{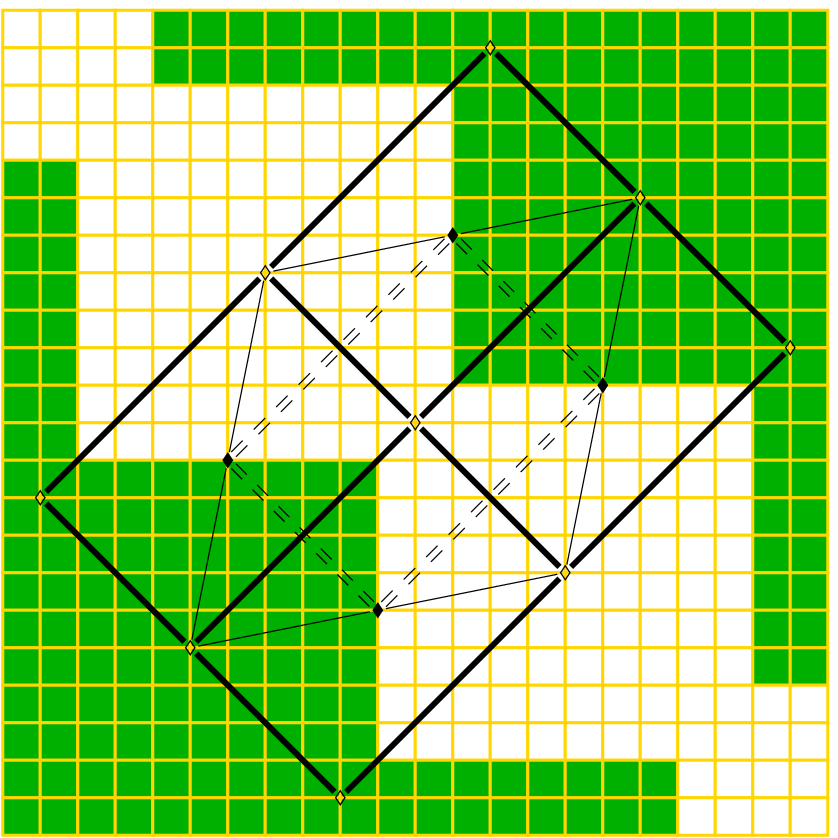}
\caption{Order-48 fabric of species 29.}
\label{fig:14a}
\end{figure}

\begin{figure}
\[
\begin{array}{cc}
\epsffile{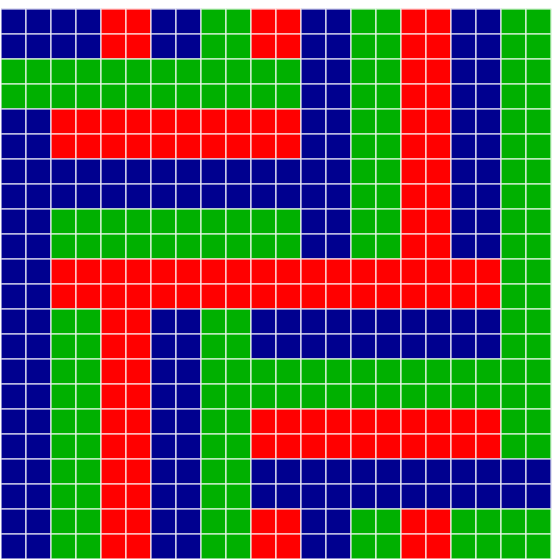} &\epsffile{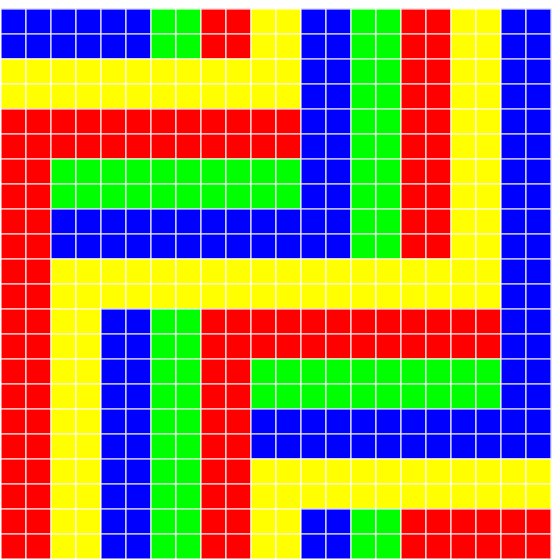}\\
\mbox{(a)} &\mbox{(b)}
\end{array}
\]
\caption{Thick stripings of fabric of Fig.~\ref{fig:14a}. a. With three colours. b. With four colours.}
\label{fig:14a_33;61b}
\end{figure}
An example from species 29 is the fabric of Fig.~\ref{fig:14a}, 3-coloured by thick striping in Fig.~\ref{fig:14a_33;61b}a.

For subspecies $27_o$ and species 30, glides in both perpendicular directions are odd.
Axes that can be an odd multiple ($\geq 1$) of $3\delta$ apart can be placed along the axes between the dark blocks, and the mirrors can be placed along the mirrors in the dark blocks.
The fabric's perpendicular axes and mirrors can then be placed along the perpendicular axes and mirrors respectively.
An example of species 30 is the fabric of Fig.~\ref{fig:13a;7a}a, 3-coloured by thick striping in Fig.~\ref{fig:13a;7a}b.

The dimensions $3\delta$ and $6\delta$ for three colours need only be changed to $p\delta$ and $2p\delta$ for any odd number of colours $p$ to specify how thick striping can be done in each possible species.
\begin{Thm}
In each of the species $11, 13, 15, 17, 19, 21, 22, 25, 27, 29$, and $30$ there are fabrics that can be perfectly coloured by thick striping with an odd number $2m+1$ of colours and redundant cells arranged as a doubled $2m/1$ twill for $m=1, 2, \dots$.
\label{thm:4}
\end{Thm}
\begin{figure}
\[
\begin{array}{cc}
\epsffile{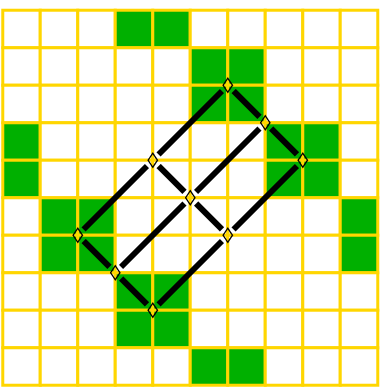} &\epsffile{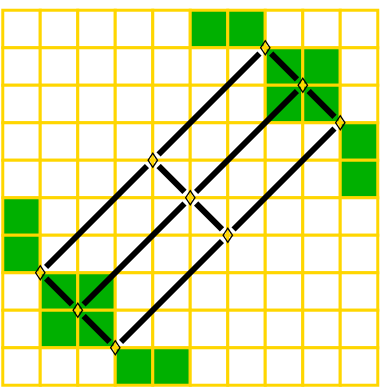}\\
\mbox{(a)} &\mbox{(b)}
\end{array}
\]
\caption{a. The 8-3-1 doubled twill as redundant blocks for 4-colouring with thick stripes. b. The 12-3-1 doubled twill as redundant blocks for 6-colouring with thick stripes.}
\label{fig:59e;59f}
\end{figure}

\section{Four and Six colours}
\label{sect:4&6}

\begin{figure}
\[
\begin{array}{cc}
\epsffile{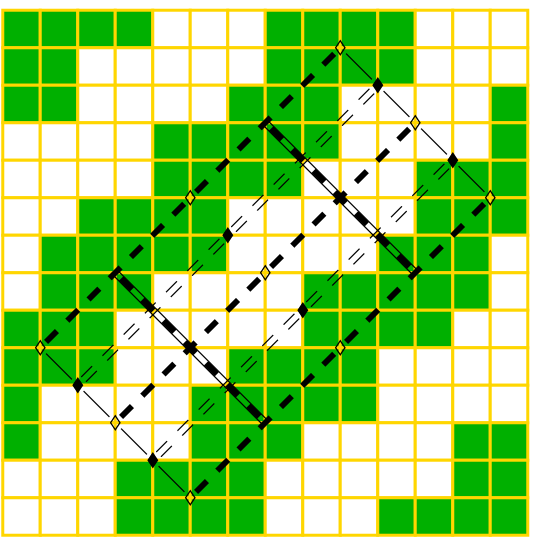} &\epsffile{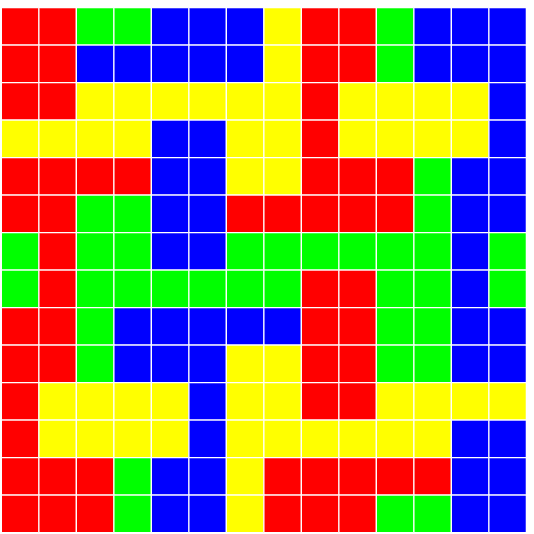}\\
\mbox{(a)} &\mbox{(b)}
\end{array}
\]
\caption{a. Fabric of species 21. b. Four-colouring by thick striping.}
\label{fig:60a;60b}
\end{figure}
\begin{figure}
\[
\begin{array}{cc}
\epsffile{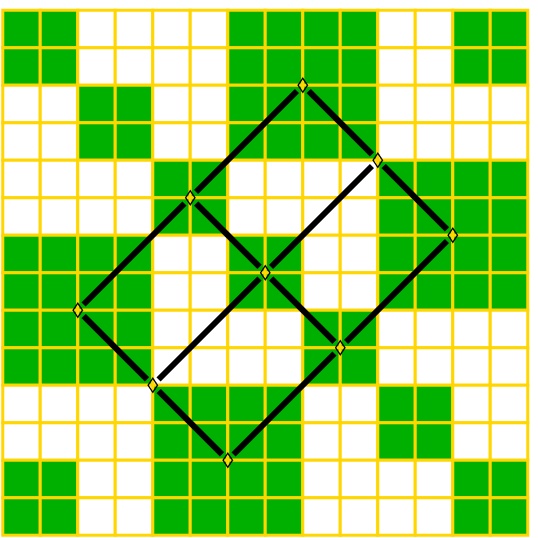} &\epsffile{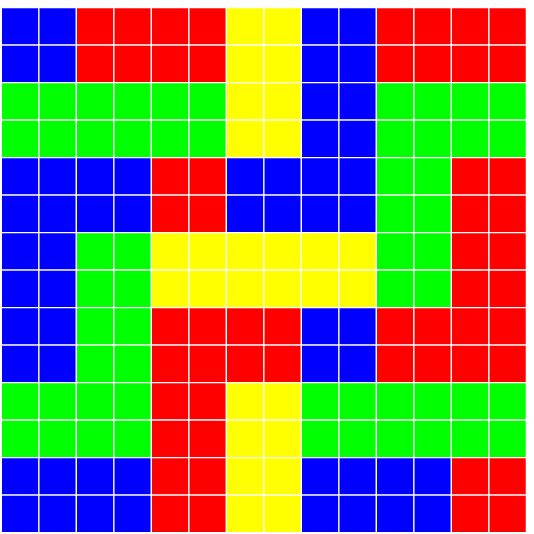}\\
\mbox{(a)} &\mbox{(b)}
\end{array}
\]
\caption{a. Subspecies-$25_e$ example. b. Four-colouring by thick striping.}
\label{fig:22a;60c}
\end{figure}

\noindent As remarked in Section \ref{sect:Finite}, the redundant doubled twill for thick striping with four colours is 8-3-1 of subspecies $25_e$ (Fig.~\ref{fig:59e;59f}a).
Mirrors are spaced $\delta$ and $2\delta$ apart in perpendicular directions.
As for three colours, the twillin and quarter-turn bans for thick striping reduce the possible species to 11, 13, 15, 17, 19, 21, 22, 25, 27, 29, and 30.
Subspecies  $11_o$, $13_o$, $15_o$, $17_o$, $19_o$, $25_o$, $27_o$, and 30 have only odd glides, and even subspecies  $11_e$, $13_e$, $15_e$, and 22 have some glides odd in $\delta$, and no such glide-reflection can be accommodated in the symmetry group of 8-3-1.
The remaining subspecies, $17_e$, $19_e$, 21, $25_e$, $27_e$, and 29 can be thickly striped with four colours.
All but $25_e$ combine glide-reflections and mirrors, and the mirrors of $25_e$ can be regarded as glide-reflections with zero glide.
All glides, which are inter-mirror distances, must be even and can be so, provided that the axes of glide-reflection are placed along mirrors perpendicular to the dark lines and mirrors are placed, say, on the dark lines.
That in turn necessitates those glides' being multiples ($\geq 0$) of $4\delta$.
Perpendicular glide-reflections (species $27_e$, and 29) with axes intervening between the fabric's mirrors, spaced out at multiples of $4\delta$ by the glide-reflections (or without them in the case of $25_e$), say in the dark lines, fall on the mirrors parallel to them between the dark lines.
Fabrics of subspecies $27_e$ and 29 have central rectangle twice an odd multiple of $\delta$ by twice an even multiple of $\delta$, and so there is only one way that their axes can be arranged too.
An example of the configuration of species $17_e$, $19_e$ and 21 is a fabric of species 21 illustrated in Fig.~\ref{fig:60a;60b}a and 4-coloured in Fig.~\ref{fig:60a;60b}b.
Examples of species $25_e$ and 29 are the fabrics in Figg.~\ref{fig:22a;60c}a and \ref{fig:14a}, 4-coloured in Figg.~\ref{fig:22a;60c}b and \ref{fig:14a_33;61b}b respectively.

As for three colours, the dimension $4\delta$ for four colours need only be changed to $4p\delta$ and for any number of colours $4p$ to specify how thick striping can be done in each possible species.

\begin{figure}
\epsffile{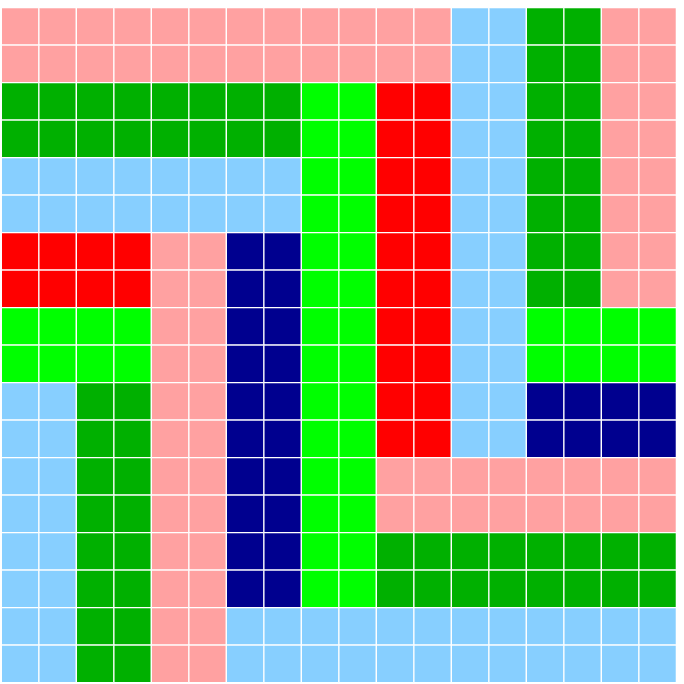}
\caption{Six-colouring by thick striping of the species-$21$ fabric of Fig.~\ref{fig:18a;62_3}a.}
\label{fig:62_6}
\end{figure}

\begin{figure}
\epsffile{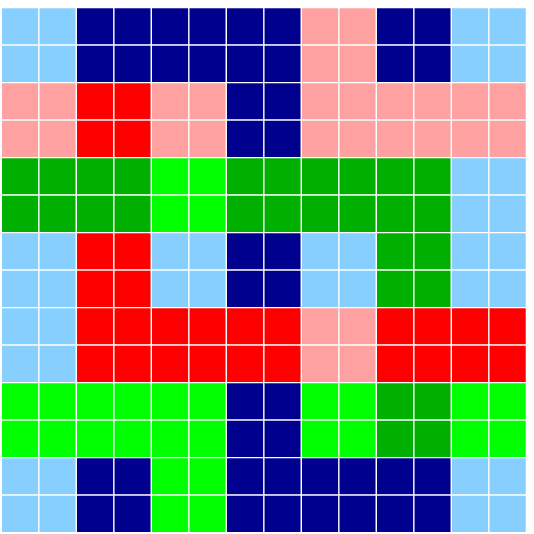}
\caption{Six-colouring by thick striping of the species-$25_e$ fabric of Fig.~\ref{fig:22a;60c}a.}
\label{fig:63}
\end{figure}

\begin{figure}
\epsffile{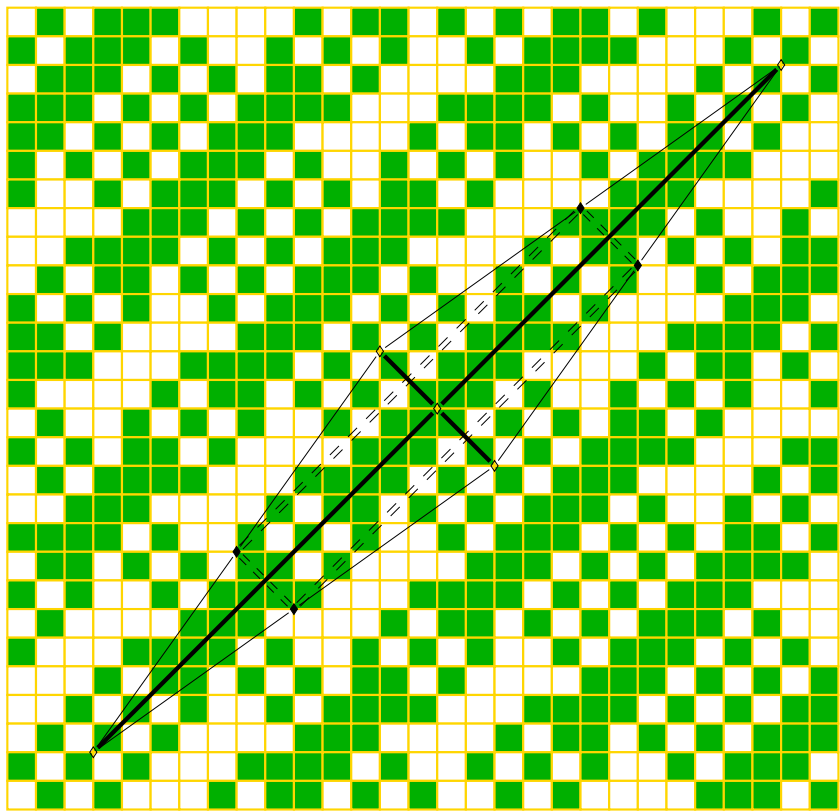}
\caption{Species-$29$ fabric.}
\label{fig:64}
\end{figure}

\begin{figure}
\epsffile{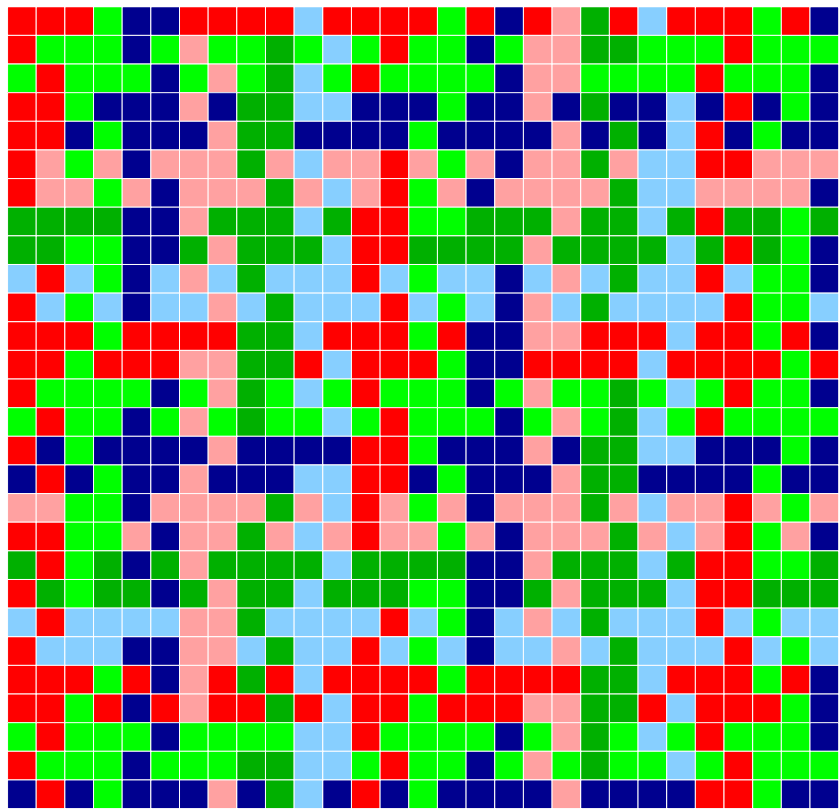}
\caption{Six-colouring by thick striping of the species-$29$ fabric of Fig.~\ref{fig:64}.}
\label{fig:65}
\end{figure}

The redundant doubled twill for thick striping with six colours is the 5/1 twill doubled, 12-3-1 of subspecies $25_o$ with basic rectangle (quarter lattice unit) $\delta$ by $3\delta$ (Fig.~\ref{fig:59e;59f}b).
The twillin and quarter-turn bans reduce the possible colorable species to 11, 13, 15, 17, 19, 21, 22, 25, $27_e$ and 29, $27_o$ and 30 as for four colours, but the oddness of both dimensions of the basic rectangle make it only very slightly easier to fit in symmetry groups.
Glide-reflections that can be imbedded in $G_2$ of 12-3-1 must still have even glides in both perpendicular directions to take redundant blocks to redundant blocks.
Accordingly, species or subspecies 11, 13, 15, $17_o$, $19_o$, 22, $27_o$, and 30 are eliminated as all require at least one glide to be odd.
The only difference from the four-colour list ($17_e$, $19_e$, 21, $25_e$, $27_e$, 29) turns out to be the species $25_o$ of the redundant twill itself.
No species is added to the list, just that subspecies.
All of the possibilities can be coloured.
Consecutive axes of glide-reflection ($17_e$, $19_e$, 21) that can be an odd multiple ($\geq 1$) of $3\delta$ apart can be placed on and between dark lines' mirrors, their perpendicular mirrors, an even multiple of $\delta$ apart falling on mirrors perpendicular to the dark lines.
An example of species 21 is the fabric of Fig.~\ref{fig:18a;62_3}a, illustrated 6-coloured by thick striping in Fig.~\ref{fig:62_6}.
Species 25 has only mirrors, and so consecutive mirrors, if an odd multiple of $3\delta$ apart, can be placed on mirrors along and between the dark lines, the perpendicular mirrors falling on mirrors perpendicular to the dark lines regardless of subspecies.
An example of subspecies $25_e$ is the fabric of Fig.~\ref{fig:22a;60c}a, illustrated 6-coloured by thick stiping in Fig.~\ref{fig:63}.
The glide-reflection axes ($17_e$, $19_e$, 21) could alternatively be placed along mirrors perpendicular to the dark lines.
That requires the glides' being multiples ($\geq 1$) of $6\delta$.
Also requiring these glides are species $27_e$ and 29, where mirrors a multiple of $6\delta$ apart can be placed along the dark lines.
Glide-reflections with axes between those mirrors then fall on the mirrors between the dark lines, and the glide-reflection axes and mirrors perpendicular to them lie on perpendicular mirrors.
An example of a fabric of species 29 is illustrated in Fig.~\ref{fig:64} and 6-coloured in Fig.~\ref{fig:65}.
The example is more interesting than attractive.
Its comparative unattractiveness results not just from the large number of colours but also from the fact that so much of the symmetry of the fabric is lost in the colouring. 
The translations along the edges of the lattice unit are side-reversing and so disappear in the coloured fabric along with the side-reversing half-turns and the mirrors.

For $4p+2$ colours, $3\delta$ and $6\delta$ above need to be replaced by $(2p+1)\delta$ and $(4p+2)\delta$.
\begin{Thm}
In each of the species $17, 19, 21, 25, 27$, and $29$ there are fabrics that can be perfectly coloured by thick striping with  $2m$ colours and redundant cells arranged as a doubled $(2m-1)/1$ twill for $m= 2,3, \dots$.
\label{thm:5}
\end{Thm}

\section{Colouring woven Tori}
\label{sect:tori}

All of the plane fabrics whose colourings are considered in sections \ref{sect:three colours} and \ref{sect:4&6} can be used to weave flat tori as discussed at the end of section \ref{sect:Background}. For the tori to be coloured, the topological constraint that the outline of the torus must be a period parallelogram of the weaving remains in place, but there is also the colouring constraint that the torus must be a period parallelogram of the planar colouring. 

Given a perfect $c$-colouring by thick striping of an order-$n$ isonemal design, one can consider how the plane can be mapped to a torus, that is, how a suitable period parallelogram can be chosen. Since the design is given, with its various period parallelograms, both $n$-by-$n$ and oblique $H_1$ lattice units, the question is how many of each must be combined to form a period parallelogram of the pattern. In the  $n$-by-$n$ case, the smallest $mn$-by-$mn$ choice must have $mn$ be the lowest common multiple of $n$ and $2c$, since each stripe of the colouring is two cells wide. The size of a lattice unit must be inflated similarly. The smallest period parallelograms of any colouring can of course be assembled to make others as large as one pleases.

Another approach would take the number of colours $c$ and the species as given and select the lattice unit appropriate to the species, if any, so that fabrics based on it are perfectly colourable.

Consider the example of Fig.~\ref{fig:59c;59d}b, a 3-coloured species-$19_e$ fabric of order 24. The lattice unit shown in Fig.~\ref{fig:59c;59d}a is too small to become a torus with the colouring of Fig.~\ref{fig:59c;59d}b because it has five stripes in each direction. The rectangle referred to in \cite{Thomas2010a} as basic for species $19_e$ --- either the smallest outlined by glide-reflection axes or with centres of half-turns at corners --- is $2\delta$ by $3\delta$. It is one quarter of the period parallelogram that is the $H_1$ lattice unit. In general, the number of stripes is $a+b$ where the basic rectangle is $a\delta$ by $b\delta$. It requires three of these lattice units in the direction of the even dimension of the basic rectangle to make the colours come out even. The $H_1$ lattice unit mapping to the torus, $4\delta$ by $6\delta$, must be tripled in the SW-NE direction to $12\delta$ by $6\delta$. In this case, the $mn$-by-$mn$ square of the second paragraph back would be 24 by 24 since $m$ can be 1, with area 576 as against $12\delta \times 6\delta =144$ for the oblique rectangle. More interestingly, the $24\times 24$ square/torus has four thick stripes (separate strands) of each colour in each direction, whereas the rectangle/torus has only a single stripe of each colour in each direction, each crossing the torus three times on account of the way their ends in the rectangle match up in the torus.

\nocite{*}
\bibliographystyle{amsplain}

\begin{thebibliography}{22}


\bibitem[1]{Brahana} Brahana, H.R., `Regular maps on an anchor ring', \emph{Amer.~Math.~J.} \textbf{48} (1926), 225--240.

\bibitem[2]{Coxeter} Coxeter, H.S.M., `Configurations and maps', \emph{Reports of a Mathematical Colloquium} (2) \textbf{8} (1948), 18--38.

\bibitem[3]{Emery} Emery, Irene, \emph{The primary structures of fabrics}. Second revised edition. Washington, D.C.: The Textile Museum, 1980. First edition 1966, Thames \& Hudson reprint 1994.

\bibitem[4]{GS1980} Gr\"unbaum, Branko, and Geoffrey C.~Shephard, `Satins and twills: An introduction to the geometry of fabrics', \emph{Mathematics Magazine} \textbf{53} (1980), 139--161.

\bibitem[5]{GS1985} Gr\"unbaum, Branko, and Geoffrey C.~Shephard, `A catalogue of isonemal fabrics', in \emph{Discrete Geometry and Convexity}, Jacob E.~Goodman \emph{et al.,} eds. \emph{Annals of the New York Academy of Sciences} \textbf{440} (1985), 279--298.%

\bibitem[6]{GS1986} Gr\"unbaum, Branko, and Geoffrey C.~Shephard, `An extension to the catalogue of isonemal fabrics', \emph{Discrete Mathematics} \textbf{60} (1986), 155--192.

\bibitem[7]{GS1988} Gr\"unbaum, Branko, and Geoffrey C.~Shephard, `Isonemal fabrics', \emph{American Mathematical Monthly} \textbf{95} (1988), 5--30.%

\bibitem[8]{HT1991} Hoskins, J.A., and R.S.D.~Thomas, `The patterns of the isonemal two-colour two-way two-fold fabrics', \emph{Bull.~Austral.~Math.~Soc.} \textbf{44} (1991), 33--43.

\bibitem[9]{Roth1993} Roth, Richard L., `The symmetry groups of periodic isonemal fabrics', \emph{Geometriae Dedicata} \textbf{48} (1993), 191--210.


\bibitem[10]{Roth1995} Roth, Richard L., `Perfect colorings of isonemal fabrics using two colors', \emph{Geometriae Dedicata} \textbf{56} (1995), 307--326.

\bibitem[11]{Schatt1978} Schattschneider, D., `The plane symmetry groups: Their recognition and notation', \emph{American Mathematical Monthly} \textbf{85} (1978), 439--450.

\bibitem[12]{Thomas2009} Thomas, R.S.D., `Isonemal Prefabrics with Only Parallel Axes of Symmetry', \emph{Discrete Mathematics} \textbf{309} (2009), 2696--2711. Online: \url{http://arxiv.org/abs/math/0612808v2} and on the journal site, doi:10.1016/j.disc.2008.06.028.

\bibitem[13]{Thomas2010a} Thomas, R.S.D., `Isonemal Prefabrics with Perpendicular Axes of Symmetry', \emph{Utilitas Mathematica} \textbf{82} (2010), 33--70. Online: \url{http://arxiv.org/abs/0805.3791v1}.

\bibitem[14]{Thomas2010b} Thomas, R.S.D., `Isonemal Prefabrics with no Axes of Symmetry', \emph{Discrete Mathematics} \textbf{310} (2010), 1307--1324. doi:10.1016/j.disc.2009.12.015. Online: \url{http://arxiv.org/abs/0911.1467v2} and on the journal site.

\bibitem[15]{Thomas2011} Thomas, R.S.D., `Perfect colourings of isonemal fabrics by thin striping', \emph{Bull.~Australian Math.~Soc.} \textbf{83} (2011), 63--86. 
Online: \url{http://arxiv.org/abs/1006.5653} and on the journal site, doi:10.1017/S0004972710001632. 

\bibitem[16]{Thomas2012} Thomas, R.S.D., `Perfect colourings of isonemal fabrics by thick striping', \emph{Bull.~Australian Math.~Soc.} \textbf{85} (2012), 325--349. 
Online: \url{http://arxiv.org/abs/1109.2254} and on the journal site, doi:10.1017/S0004972711002899. 
\end{thebibliography}

\end{document}